# DISCUSSION: CONDITIONAL GROWTH CHARTS

By Raymond J. Carroll[1] and David Ruppert[2]

*Texas A&M University and Cornell University*

**1. Overview.** Wei and He are to be congratulated on an innovative and important article. The conditional approach to growth charts described in their article is important in a practical sense, and the use of quantile regression is both natural and well motivated. We look forward to further application of their idea to actual practice, because the concept of "falling behind" in one's growth cycle has two meanings: the usual standard growth chart, and the conditional growth chart described here. We have described the Wei and He approach to pediatricians, and they all grasped the essential clever idea immediately and were enthusiastic about the idea.

Our commentary will focus on three aspects of the approach used by the authors, specifically (a) the use of *unpenalized* B-splines as described by the authors; (b) conditional versus marginal semiparametric modeling of longitudinal data; and (c) some alternative modeling approaches to "catch-up" that may get at the issue more directly and flexibly.

**2. B-splines should be penalized.** One purpose of discussions, of course, is to make things a little lively, and here is our contribution. In our view, one should have some skepticism of how nonparametric unpenalized B-splines and unpenalized regression splines really are in the context of *non*parametric regression. More precisely, and less inflammatory, the connection between asymptotic theory for unpenalized splines and actually attempting to be at least reasonably nonparametric is not at all clear.

There is obviously a need to balance practical behavior and ease of implementation with theory. Kernel methods (see below) are one means of doing this. In the spline context, there are four approaches: smoothing splines,

Received November 2005.
[1]Supported by National Cancer Institute Grant CA-57030 and by the Texas A&M Center for Environmental and Rural Health via National Institute of Environmental Health Sciences Grant P30-ES09106.
[2]Supported by NSF Grant DMS-04-34390 and NIH Grant CA-57030.

penalized regression splines, unpenalized regression splines and free-knot splines; see [18] for a recent review. As authors, we have observed the following.

- Smoothing splines are basically penalized regression splines that place a knot at every value of the covariate. We have no idea how to do asymptotic theory for smoothing splines in the Wei and He context, and we would be interested to see if such a theory is even possible, for example, along the lines of [13].
- Lower-order knots penalized regression splines are rapidly becoming a practical method of choice. However, every time we write a paper about this technique, we have been reviewed by a smoothing spline person who asks the questions: "how many knots" and "where do you place the knots." In response, Ruppert [16] did an extensive numerical study and showed that, in effect, penalized regression splines with data-adaptive penalties and 40 knots will work quite well in many practical settings. Ruppert and Carroll [17] and Ruppert, Wand and Carroll [18], Chapters 5 and 17, describe strategies for selecting the number of knots. There is little in the way of asymptotic theory for penalized regression splines, unfortunately, presumably for the same reason that penalized smoothing splines are difficult to analyze; this point is also brought up in every review. For some theory, see [3].
- Unpenalized regression splines with 40 knots are generally ridiculously non-smooth, which is of course why penalties are used. For example, in Figure 1, we generated data according to the regression model $Y = \sin(2X) + \varepsilon$, where $X$ is equally spaced on the unit interval, $n = 400$ and $\varepsilon = \text{Normal}(0,1)$. We used 40 knots and fit penalized and unpenalized cubic B-splines, using the Matlab software available from Brian Marx at www.stat.lsu.edu/faculty/marx/.
- There is a substantial literature on unpenalized B-splines that achieve their smoothness via knot selection, the so-called free-knot spline methodology. Much of this literature is Bayesian, being based on model averaging; see, for example, [2, 4, 19]. With such complex methods, theoretical results are naturally generally lacking, although see [12] and [20] for recent non-Bayesian approaches with impressive (asymptotic) theoretical treatment.

In summary, the nonparametric spline literature generally uses a fairly large number of knots, realizes that penalties must be imposed in one way or another to obtain smoothness, and achieves this smoothness either by direct penalties, or by some version of knot selection. Often for these methods asymptotic theory is not available.

In contrast, unpenalized regression spline methods do have a beautiful asymptotic theory, as exemplified here and in a series of important papers by Jianhua Huang and colleagues [5, 6, 7, 8, 9]. Plots, such as in Wei and He, are often pleasantly smooth, in contrast to Figure 1. The obvious question is: what is going on?

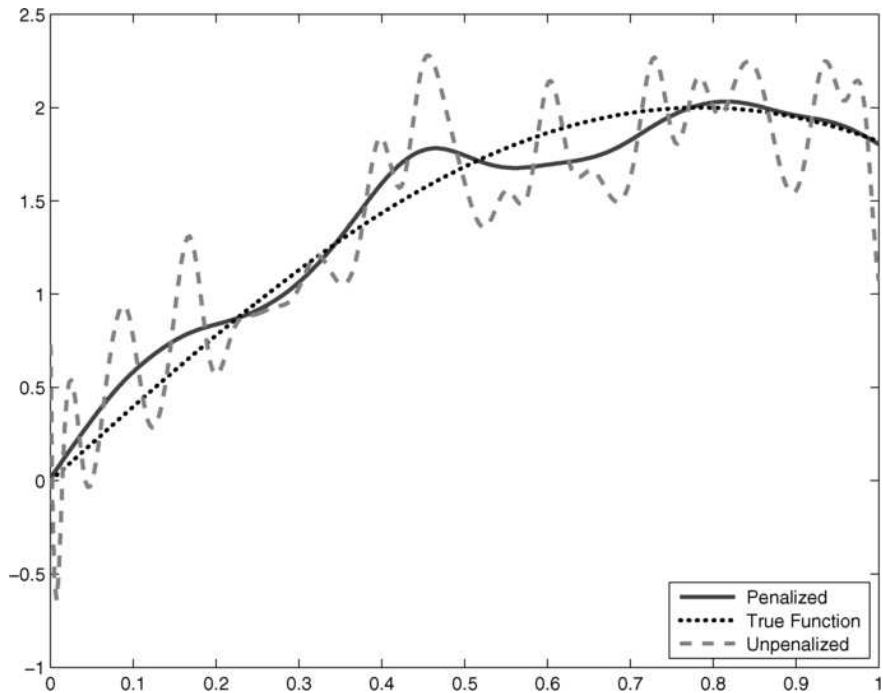

FIG. 1. *Fits to the regression model $Y = \sin(2X) + \varepsilon$, where $X$ is equally spaced on the unit interval and $\varepsilon = \text{Normal}(0,1)$. The dotted line is the true function. The solid line is a penalized cubic B-spline fit with* 40 *knots. The dashed line uses the same basis, but lacks a penalty.*

To get such smoothness in an unpenalized regression spline, one necessarily must insist that the number of knots be small. Indeed, Wei and He use cubic B-splines with three knots in their examples. This is no fluke: if one fixes the knots, and does no penalization or knot selection, there is no getting around phenomena such as in Figure 1. Indeed, the basic point is actually one of the conditions of their Theorem 3.1. Note that the essential condition for a function of bounded second derivative is that the number of knots be proportional to $k \propto n^{1/5}$, or equivalently, that the sample size $n \propto k^5$. Since $3^5 = 243$, $4^5 = 1024$, $5^5 = 3125$ and $6^5 = 7776$, we see that the asymptotics essentially require that for any practical problem, the number of knots should be no more than six. This is hardly nonparametric regression! Of course, these calculations are deliberately shocking and totally slanted in order to add some controversy, because we have not mentioned constants of proportionality, but once one does get into estimating the number of knots and where they should be placed, then it is not clear what the theory would be.

We note that Wei and He's criterion function (2.2) could have been analyzed by kernel methods, using local-likelihood ideas. Theorem 3.1 would then have been easy to analyze using either profiling or backfitting; see, for example, [1] for the computationally far easier backfitting and [21] for the more complex profiling. We conjecture, of course, that the same limiting result as in Theorem 3.1 would have been obtained, since one way to interpret the work of Newey [15] is that different implementations of the same criterion function should not lead to different limiting results for parametric components. It would be interesting to know whether this conjecture is correct.

**3. Marginal versus conditional longitudinal models.** Wei and He deal with longitudinal semiparametric models, using a conditional approach that is ideally suited to their main idea.

The usual approach to this problem is via a marginal model, for example, in a slightly altered notation,

$$Y_{ij} = m(X_{ij}, \beta) + \theta(Z_{ij}) + \varepsilon_{ij};$$
(1)
$$\text{cov}(\varepsilon_{i1}\ldots) = \Sigma_i(\tau),$$

where $m(\cdot)$ is a known function. Wang, Carroll and Lin [22] describe the semiparametric efficient solution to this problem if the $\varepsilon_{ij}$ are Gaussian using kernels, although the work of Lin et al. [11] essentially shows the same result for smoothing splines. There is considerable controversy about this problem. Early solutions (and many later ones) simply ignored the correlation structure, with a notable exception being Zeger and Diggle [23]. The semiparametric efficient solution for Gaussian data is explicit, that is, not iterative, but it still takes some work to implement.

Wei and He's conditional model approach, in contrast, neatly avoids all these issues, because their model is

$$(2) \qquad Y_{ij} = m(Y_{i,k}, X_{ij}, \beta^*) + \theta(Z_{ij}) + \varepsilon_{ij}^*,$$

where $k < j$ and the $\varepsilon_{ij}^*$ are independent. Here even the kernel approach is simple to implement.

We would be interested in Wei and He's thoughts on how to use model (2) to help understand marginal models such as (1), *in the context of their innovative quantile regression modeling.* Since they (quite properly) do not work in a mean-based model, the marginal interpretation of their model, if any, is unclear to us. If there were a marginal interpretation, then perhaps this could be used in the more standard unconditional growth chart arena.

On a theoretical note, if one starts from a conditional model and then turns it into a marginal model, it is not clear whether the methods of estimation are semiparametric efficient *in the marginal model*. Recent work by Lin and Carroll [10] can be used to answer this question.

**4. Other literature.** Much work on nonparametric estimation of marginal curves has been done in the past using kernel methods; see [14] and the references therein. We are uncertain how extensively these methods have been adopted by practitioners, but since they are not mentioned in this paper it seems that they are at least somewhat neglected. Perhaps Professors Wei and He can comment. In Müller's book there is an emphasis on estimating the first two derivatives, that is, growth velocity and acceleration. Müller mentions growth spurts during adolescence and perhaps these might be missed by a three-knot cubic spline. Therefore, we wonder whether the methods in this paper can be applied to long time spans that include adolescence or, instead, whether a methodology with a smoothing parameter would be needed.

In (5.1), it seems more appropriate to use $H_{i,j}^3$ instead of $H_{i,j}$ as a covariate, since weight should be roughly proportional to the cube of a linear dimension. Perhaps the increase in $c_\tau$ with $\tau$ seen in Table 2 and mentioned by the authors is due to using $H_{i,j}$ instead of its cube. Similarly, the authors discuss a boy who jumped from the 0.25th quantile to the median level at age 0.61 year and concluded "that he might be overweight given his prior path and current height." This is an interesting and potentially important conclusion, but one must be certain that the effect of height is being modeled correctly before making it.

**5. Alternative models with catch-up.** In their Introduction, Professors Wei and He mention "catch-up" growth where subjects on the upper (or lower) centiles move toward the median at a faster rate than others. Model (2.1) is very general and should include the possibility of catch-up growth, but it would be helpful to practitioners if the mechanism in the model for catch-up growth was explicit. Consider (5.1) of the paper with $k = 1$, that is,

$$W_{i,j} = g_\tau(t_{i,j}) + (a_\tau + b_\tau D_{i,j})W_{i,j-1} + c_\tau H_{i,j} + e_{i,j}.$$

There is no apparent mechanism for catch-up growth here, since, with $a_\tau$ and $b_\tau$ both positive as in Table 2, $W_{i,j}$ is an increasing function of $W_{i,j-1}$. Of course, it might be that $H_{i,j}$ catches up and forces $W_{i,j}$ to catch up as well. This raises the question of whether $H_{i,j}$ should be a covariate or a second response. Perhaps height and BMI should be modeled jointly as a bivariate dynamic process.

We have thought about other models where catch-up growth might be more explicit. One model, which we realize may be too simple but could be a good starting point, is

$$W_{i,j} = W_{i,j-1} + \{g_\tau(t_{i,j}) - g_\tau(t_{i,j-1})\}$$
(3) $$\qquad + bD_{i,j}\{W_{i,j-1} - g_\tau(t_{i,j-1})\} + D_{i,j}e_{i,j}$$
$$= W_{i,j-1} + \text{average change in the population} + \text{subject-specific change}.$$

If $b < 0$, then there is catch-up growth, because subjects with $W_{i,j-1} - g_\tau(t_{i,j-1})$ positive (negative) tend to grow less (more) than average. One feature of (3) is that $W_{i,j} \to W_{i,j-1}$ as $t_{i,j} \to t_{i,j-1}$ (so that $D_{i,j} \to 0$). Anything else, of course, would not be realistic. Besides the parameters in $g_\tau$, (3) has only a single parameter $b$. However, $b$ probably should depend on time since Professors Wei and He mention that catch-up tends to be time-specific. One might replace $b$ by $b\{(t_{i,j} + t_{i,j-1})/2\}$ where $b(\cdot)$ is a spline.

Let $W_{i,j}^* = W_{i,j} - g_\tau(t_{i,j})$. Then (3) with $b$ a spline can be written as

$$W_{i,j}^* = [1 - b\{(t_{i,j} + t_{i,j-1})/2\}D_{i,j}]W_{i,j-1}^* + D_{i,j}e_{i,j},$$

which is similar to an AR(1) process.

**6. Concluding comments.** Notwithstanding our comments about general statistical methodology, and our comments about practical details, Wei and He have written an important paper, and we look forward to reading about further developments of their innovative ideas. We have rarely read such a thought-provoking paper that has the potential to be extremely important in practice.

DEPARTMENT OF STATISTICS
TEXAS A&M UNIVERSITY
COLLEGE STATION, TEXAS 77843-3143
USA
E-MAIL: carroll@stat.tamu.edu

SCHOOL OF OPERATIONS RESEARCH
AND INDUSTRIAL ENGINEERING
225 RHODES HALL
CORNELL UNIVERSITY
ITHACA, NEW YORK 14853
USA
E-MAIL: dr24@cornell.edu